\documentclass[a4paper,12pt,reqno]{amsart}
\usepackage{amsmath}
\usepackage{amsthm}
\usepackage{amssymb}
\usepackage{latexsym}               
\usepackage{amsgen}
\usepackage[mathscr]{eucal}
\numberwithin{equation}{section}
\usepackage[dvipdfm, usenames]{color}

\newcommand{\ti}{\tilde}

\makeatother

\newtheorem{lem}{Lemma}[section]
\newtheorem{thm}[lem]{Theorem}
\newtheorem{pro}[lem]{Proposition}

\newtheorem{remark}{Remark}

\begin{document}

\title
{Existence of radially symmetric\\
stationary solutions 
for \\ viscous and Heat-conductive ideal Gas }
\author{Itsuko Hashimoto$\dagger$}
\email{itsuko@se.kanazawa-u.ac.jp}
\address{Kanazawa University, Osaka City University, Japan}

\author{Akitaka Matsumura}
\email{akitaka@math.sci.osaka-u.ac.jp}
\address{Osaka University, Japan}

\thanks{
$\dagger$ \ Supported by JSPS Grant No. 21K03306
}

\date{}
\maketitle

\def\cal#1{{\fam2#1}}


\begin{center}
\small{\it Dedicated to Professor Shih-Hsien Yu on the occasion of his 60th birthday}
\end{center}


\begin{abstract}
  
We consider the existence of radially symmetric stationary solutions of the compressible viscous and heat-conductive polytropic ideal fluid on the unbounded exterior domain of a sphere in $\mathbb{R}^n (n\ge 3)$ where the boundary and far-field conditions are prescribed.
The unique existence of the stationary solution is shown for both inflow and outflow problems in a suitably small neighborhood of the far-field state.
Estimates of the algebraic decay rate toward the far field state are also obtained. 


\end{abstract}

\medskip


\section{Introduction and Main Theorem}
\noindent  
In the present paper, we consider a Newtonian polytropic ideal model of a compressible viscous and heat-conductive fluid in the exterior domain $\Omega=\{x \in\mathbb{R}^{n}\,(n\ge 3); |x| > r_0\}$\,($r_0$ is a positive constant and $n$ is the space dimension) described by 
\begin{eqnarray}
 \label{cns} 
 \left\{
\begin{array}{l}
\rho_t+\nabla\cdot (\rho U)=0, \\[10pt]
(\rho U)_t+\nabla\cdot(\rho U\otimes U)+\nabla p =\nabla \cdot D,\quad
t>0,\ x \in \Omega,\\[10pt]
\left\{\rho \left(e + \frac{1}{2}|U|^2\right)\right\}_t +
\nabla\cdot \left\{\rho U \left(e + \frac{1}{2}|U|^2\right)\right\}+\nabla\cdot(pU)= \nabla\cdot (DU -q),\\[10pt]
\end{array}
 \right.\, 
\end{eqnarray} 
in the Eulerian coordinates, where $\rho=\rho(t,x)>0$ is the mass density, $U={}^t(u_1(t,x),\cdots ,u_n(t,x))$ is the fluid velocity, $e:=c_V \theta$, and $p:=R\rho\theta$ are internal energy per unit mass and the pressure, respectively, $\theta=\theta(t,x)>0$ is the absolute temperature. In addition, $c_V$ and $R$ are the specific heat at constant volume and the gas constant, which are assumed to be positive constants. Furthermore, $D$, $\Theta$, and $q$ are the stress tensor {\it w.r.t.} viscosity, the strain tensor, and the heat flow vector, respectively, defined by 
\begin{align}
D= 2\nu\,\Theta+\lambda\,(\nabla\cdot U)I,\quad \Theta= \frac{1}{2}\{\nabla\otimes U + {}^t(\nabla\otimes U)\},\quad q=-\kappa \nabla\theta,
\end{align}
where $\nu$ and $\lambda$ are the shear and second viscosity coefficients, and $\kappa$ is the thermal conductivity, which are assumed to be constants satisfying 
\begin{align}
\label{visc}
\nu>0, \quad 2\nu+ n\lambda \ge 0, \quad \kappa>0.
\end{align}
\noindent
In this paper, we focus our attention on the radially symmetric solutions, which have the form
\begin{align}
\label{t}
\rho(t,x)=\rho(t,r), \quad U(t,x)=\frac{x}{r}\,u(t, r),\quad \theta(t,x)=\theta(t,r), \quad r= |x|,
\end{align}
where $u(t, r)$ is a scalar function. By plugging  (\ref{t}) into (\ref{cns}), 
we can rewrite (\ref{cns}) and obtain the equations for the radially symmetric solution $(\rho, u,  \theta)$  in the form
\begin{equation}
\label{nr}
\left\{\begin{array}{l}
\rho_t + \dfrac{ (r^{n-1}\rho u)_r}{r^{n-1}}=0, \\[10pt]
(\rho u)_t
+ \dfrac{ (r^{n-1}\rho u^2)_r}{r^{n-1}}+ p_r =\mu \Big(\dfrac{(r^{n-1}u)_r}{r^{n-1}}\Big)_r, \\[10pt]
\left\{\rho \left(c_V\theta + \frac{1}{2}u^2\right)\right\}_t +
\dfrac{\left\{r^{n-1}\rho u \left(c_V\theta + \frac{1}{2}u^2\right)\right\}_r}{r^{n-1}}
+\dfrac{ (r^{n-1}pu)_r}{r^{n-1}}
\\
\hspace{3cm}= \kappa\dfrac{(r^{n-1}\theta_r)_r}{r^{n-1}} +
\mu \Big(\dfrac{(r^{n-1}u)_r}{r^{n-1}}\Big)_ru +\Psi, \quad \ \ t>0,\ r>r_0,
\\[10pt]
\Psi := 2\nu u_r^2 +2\nu(n-1)(\dfrac{u}{r})^2+\lambda(u_r+\dfrac{n-1}{r}u)^2,
\end{array}
 \right.\,
\end{equation}
where 
\begin{align}
\label{viscondi}
\mu := 2\nu + \lambda>0,
\end{align} and $\Psi$ is the dissipation function, a non-negative quadratic form  under the condition (\ref{visc}), representing the rise in temperature due to viscosity. By the definition (\ref{viscondi}), we also see that 
\begin{align}
\label{mutonu}
\frac{\nu}{\mu}\le \frac{n}{2(n-1)} \le 1.
\end{align}

Now, we consider the initial boundary value problems 
for (\ref{nr}) under the initial condition
\begin{equation}
(\rho, u, \theta)(0,r)=(\rho_0, u_0, \theta_0)(r),\quad r >r_0, 
\end{equation}
the far field condition
\begin{equation}
\displaystyle{\lim_{r \to \infty}}(\rho, u, \theta)(t,r)=(\rho_+,u_+,\theta_+),\quad t>0,
\end{equation}
and also the following two types of boundary conditions, depending on the sign of the velocity
on the boundary
\begin{equation}
\label{bc}
 \left\{
\begin{array}{ll}
(\rho, u, \theta)(t,r_0) =(\rho_-, u_-, \theta_-),\quad t>0,& (u_->0),\\[10pt]
(u, \theta)(t,r_0) =(u_-, \theta_-),\quad t>0,& (u_-\le 0),
\end{array}
 \right.
\end{equation}
where $\rho_{\pm}>0, u_{\pm}$ and $\theta_{\pm}$ are given constants. 
The case $u_->0$ is known as ``inflow problem'', the case $u_-=0$ as 
``impermeable wall problem'', and the case $u_-<0$ as ``outflow problem''. 

Regarding the impermeable wall problem, Matsumura and Nishida \cite{MNT} first considered the original exterior problem (\ref{cns}) in $\mathbb{R}^3$, where $\Omega$ is a more general exterior domain, which includes the radially symmetric case (\ref{nr}). That is, under the boundary and far field conditions $U|_{\partial \Omega}=U|_{\infty}=0$, \ $\theta_+=\theta_-$, they showed the existence of the global solution in time around the constant state $(\rho, U, \theta)=(\rho_+, 0, \theta_+)$ and its asymptotics toward the constant state for small initial perturbations. 
On the other hand, for the radially symmetric problem (\ref{nr}) in $\mathbb{R}^n$ ($n=2, 3$), Jiang \cite{J} showed the time global existence of the solution to (\ref{nr}) under the conditions $u_{\pm}=0$, \ $\theta_+=\theta_-$, and derived the partial results on the asymptotic stability toward the constant state $(\rho_+, 0, \theta_+)$ for large initial perturbations. Nakamura-Nishibata-Yanagi \cite{NN} and Nakamura-Nishibata \cite{NN0} extended the results by Jiang \cite{J} to the case with external potential forces in $\mathbb{R}^n$ ($n\ge 3$) and showed the unique existence of the stationary solution and its asymptotic stability. 

Regarding the inflow and outflow problems, Hashimoto-Matsumura \cite{HM1} study the existence 
of the radially symmetric stationary solution for the compressible Navier-Stokes equation describing the motion of viscous barotropic gas 
on the exterior domain of a sphere, 
in a suitably small neighborhood of the far field state. 
The asymptotic stability of the stationary solutions obtained in \cite{HM1} was considered in Hashimoto-Nishibata-Sugisaki \cite{HNS} for the outflow problem, and in Huang-Hashimoto-Nishibata \cite{HHN} for the inflow problem. 
In this paper, as the first step 
to extend the results to
the compressible viscous and heat-conductive ideal fluid
we shall show the existence of stationary solutions in a suitably small neighborhood of the far field state.

The stationary problem corresponding to the problem (\ref{nr})-(\ref{bc}) is 
written as
\begin{equation}
\label{nrs0}
\left\{\begin{array}{l}
\dfrac{ (r^{n-1}\rho u)_r}{r^{n-1}}=0, \\[10pt]
\dfrac{ (r^{n-1}\rho u^2)_r}{r^{n-1}}+ p_r =\mu \Big(\dfrac{(r^{n-1}u)_r}{r^{n-1}}\Big)_r,
\\[10pt]
\dfrac{\left\{r^{n-1}\rho u \left(e + \frac{1}{2}u^2\right)\right\}_r}{r^{n-1}}
+\dfrac{ (r^{n-1}pu)_r}{r^{n-1}}\\[10pt]
\hspace{4cm}\!\!= \kappa\dfrac{(r^{n-1}\theta_r)_r}{r^{n-1}} +
\mu \Big(\dfrac{(r^{n-1}u)_r}{r^{n-1}}\Big)_ru +\Psi,\quad r\ge r_0,
\\[10pt]
\displaystyle{\lim_{r\to \infty}} (\rho, u, \theta)(r) =(\rho_+, u_+, \theta_+), \\[2mm]
(\rho, u, \theta)(r_0)=(\rho_-, u_-, \theta_-)\ \ (u_->0),\quad (u, \theta)(r_0)=(u_-, \theta_-) \ \ (u_- \le 0),\\[5pt]
p= R\rho\theta, \quad e=c_V \theta.
\end{array}
 \right.\,
\end{equation}
By integrating the first equation in (\ref{nrs0}) in terms of $r$, we easily see that it holds
\begin{align}
\label{n1}
r^{n-1}\rho(r) u(r)=\epsilon, \qquad r\ge r_0, 
\end{align}
for some constant $\epsilon$, which implies by the boundary conditions that
\begin{equation}
\epsilon = r^{n-1}_0\rho_- u_-\quad(u_->0),\qquad 
\epsilon = r^{n-1}_0\rho (r_0) u_-\quad(u_-\le0). 
\end{equation}
Due to the formula (\ref{n1}), we see that if $n \ge 2$, 
\begin{align*}
u_+=\lim_{r\to \infty}u(r)= \lim_{r\to \infty}\frac{\epsilon}{r^{n-1}\rho_+}=0.
\end{align*}
Hence, it is necessary to assume $u_+=0$ for the existence of multi-dimensional
stationary solutions. 
We first note the case $u_-$=0. In this case, we can easily see from the first equation of (\ref{nrs0}) that $u=0$. By using $u=0$, we can derive from the second and third equations of (\ref{nrs0}) that 
$$
p_r =0,\quad \kappa (r^{n-1}\theta_r)_r =0,
$$
which imply
$$
\rho\,\theta=\rho_+\theta_+,\quad \theta=\ti\theta(r):=\theta_+
+(\theta_--\theta_+)(r/r_0)^{-(n-2)}.
$$
 Thus, in the case $u_-=0$, the trivial solution of (\ref{nrs0}) is given by $(\rho, u, \theta)=(\rho_+\theta_+/\ti\theta, 0, \ti\theta)$ for any $\rho_+$ and $\theta_{\pm}>0$. So, in what follows, we assume that $u_-\neq 0$.
 Now we are ready to state the main results in the present paper.
\begin{thm} 
\label{mt}
Let $n \ge 3$ and $u_+=0$.
Then, for any $\rho_+>0$ and $\theta_+>0$,
there exist positive constants $\epsilon_0$ and $C$ satisfying the following:

\medskip

\noindent
{\rm (I)} \ Let $u_->0$. If $|u_-|+|\rho_--\rho_+|+|\theta_--\theta_+| \le \epsilon_0$, 
there exists a unique smooth solution $(\rho,u, \theta)$ of 
the problem {\rm (\ref{nrs0})} satisfying 
\begin{align}
\label{1}
\begin{aligned}
  & |\rho(r)-\rho_+|,\ |\theta(r)-\theta_+|\le Cr^{-(n-2)}(|u_-|+|\rho_--\rho_+|+|\theta_--\theta_+|),\\[5pt]
  & C^{-1}r^{-(n-1)}|u_-| \le |u(r)| \le Cr^{-(n-1)}|u_-|, \qquad r\ge r_0.
\end{aligned}
\end{align}
{\rm (II)} \ Let $u_-< 0$. If $|u_-|+|\theta_--\theta_+|\le \epsilon_0$, 
there exists a unique smooth solution $(\rho, u, \theta)$ 
of the problem {\rm (\ref{nrs0})} satisfying
\begin{align}
\label{decayoutflow}
\begin{aligned}
  & |\rho(r)-\rho_+|,\ |\theta(r)-\theta_+|\le Cr^{-(n-2)}(|u_-|+|\theta_--\theta_+|), \\[5pt]
  & C^{-1}r^{-(n-1)}|u_-| \le |u(r)|\le Cr^{-(n-1)}|u_-|, \qquad r\ge r_0.
\end{aligned}
\end{align}
\end{thm}

\medskip

\begin{remark}
When $u_-=0$, as we stated above, the trivial solution is given by $(\rho, u, \theta)=(\rho_+\theta_+/\ti\theta, 0, \ti\theta)$.
Therefore, when $u_->0$ and $\rho_- \neq \rho_+$, since 
$\rho_+\theta_+/\ti\theta(r_0) = \rho_+\theta_+/\theta_-$ is not $\rho_-$ in general, we can expect that a boundary layer for mass density does appear as $u_- \to +0$. 
\end{remark}
\begin{remark}
As for the stationary solutions for exterior problems of the viscous barotropic  model investigated in \cite{HM1}, the stationary solution of density $\rho(r)$ has the decay rate properties:
\begin{eqnarray}
\label{decaybarot}
\left\{
\begin{array}{l}
\displaystyle{|\rho(r)-\rho_+|\le Cr^{-2(n-1)}(|u_-|+|\rho_--\rho_+|)},
\quad u_->0,\quad r>r_0,
\\[5mm]
\displaystyle{|\rho(r)-\rho_+|\le Cr^{-2(n-1)}|u_-|},
\quad u_-\le 0,\quad r>r_0.
\end{array}
 \right.\,
\end{eqnarray}
We can see that the decay rate of the results in {\rm Theorem \ref{mt}} is quite different from that 
in (\ref{decaybarot}) because of the influence from the absolute temperature. 
\end{remark}

\begin{remark}
We discussed the inviscid limit problem for the barotropic model in \cite{HM2}. So, the inviscid limit problems for the full system {\rm (\ref{nrs0})} with $\mu, \ \kappa \to 0$ seem a very interesting problem in the future. In particular, as for the inviscid 
limit problem $\mu \to 0$,
we shall see in the proofs of {\rm Theorem 1.1} that if the range of $\mu$ is bounded, say $0<\mu\le 1$, then the constants $\epsilon_0$ and $C$ in Theorem 1.1 can be chosen independent of $\mu$. 

\end{remark}

\section{preliminary} 
In this section, as stated above, we assume $u_- \ne 0$. 
In what follows, we also assume $r_0=1$ without loss of generality. 
To reformulate the problem (\ref{nrs0}), we introduce the specific volume $v$ by $v=1/\rho$ (accordingly, denote $v_\pm$ by $1/\rho_\pm$). Using the first equation of  (\ref{nrs0}), the velocity $u$ is given in terms of $v$ as
\begin{align}
\label{utov}
u(r)=\frac{\epsilon}{r^{n-1}}v(r),\quad r\ge 1,
\end{align}
where $\epsilon = u_-/v_-\ (u_->0)$, and
$\epsilon = u_-/v(1)\ (u_-< 0)$.
We further introduce 
a new unknown functions $\eta$ and $\chi$, as the deviation of $v$ and $\theta$ from the far field state $v_+$ and $\theta_+$, respectively, by
\begin{equation}
\label{transec}
\eta(r) = v(r) - v_+,\qquad \chi=\theta(r)-\theta_+, \qquad r \ge 1. 
\end{equation}
Substituting (\ref{utov}) into the second equation of (\ref{nrs0}) by using (\ref{transec}), and integrating the resultant equation with respect to $r$ on $[r, \infty)$, we can rewrite the second equation of (\ref{nrs0}) as 
\begin{equation}
\label{etar}
\begin{array}{l}
\eta_r+\dfrac{R\theta_+r^{n-1}}{\epsilon \mu v_+^2}\eta = \dfrac{Rr^{n-1}}{\epsilon  \mu v_+}\chi +F(\eta,\chi),
\end{array}
\end{equation}
where, 
\begin{equation}
\label{etaF}
\begin{array}{l}
F(\eta,\chi):=\dfrac{Rr^{n-1}\chi\eta}{\epsilon  \mu v_+(v_++\eta)}
+\dfrac{R\theta_+r^{n-1}\eta^2}{\epsilon \mu v_+^2(v_++\eta)}
+\dfrac{\epsilon v_+}{2\mu r^{n-1}}+\dfrac{\epsilon\eta}{\mu r^{n-1}}\\[15pt]
\hspace{4cm}-\dfrac{\epsilon(n-1)r^{n-1}}{\mu}\displaystyle{\int_r^{\infty}\dfrac{\eta(s)}{s^{2n-1}}ds}.
\end{array}
\end{equation}
Here, we used the far field conditions, and extracted the linear parts at the far field state, on the same line as in \cite{HM1}.
To rewrite the third equation of (\ref{nrs0}), we multiply $r^{n-1}$ by it and substitute $u=\epsilon(v_++\eta)/r^{n-1}$ into the resultant equation. Then we obtain the following equation: 
\footnotesize
\begin{align}
\left\{\epsilon \left(c_V\theta + \frac{1}{2}\left(\frac{\epsilon (v_++\eta)}{r^{n-1}}\right)^2\right)\right\}_r
+(R\epsilon\theta)_r
= \kappa(r^{n-1}\theta_r)_r+\epsilon^2\mu \Big(\dfrac{(v_++\eta)\eta_r}{r^{n-1}}\Big)_r
+\Phi(\eta,\eta_r).
\end{align}
\normalsize
By integrating the above equation with respect to $r$ on $[r, \infty)$, we have
\begin{align}
\epsilon c_P\chi + \dfrac{\epsilon^3}{2}\Big(\dfrac{v_++\eta}{r^{n-1}}\Big)^2
- \kappa r^{n-1}\chi_r = \alpha
+\epsilon^2\mu \dfrac{(v_++\eta)\eta_r}{r^{n-1}}
-\int^\infty_r\Phi(\eta,\eta_r)(s) \,ds,
\end{align}
where $c_P:=R+c_V$ is the specific heat at constant pressure, $\alpha$ is a constant of integration, and
\begin{align}
\label{Phi}
\Phi(\eta,\eta_r)(r):=-4\epsilon^2\nu (n-1)\dfrac{(v_++\eta)\eta_r}{r^{n}}+2\epsilon^2\nu n(n-1)\dfrac{(v_++\eta)^2}{r^{n+1}}.
\end{align}
Thus, we finally have the following reformulated problem of (\ref{nrs0}) in terms of $(\eta, \chi)$ as
\begin{align}
\label{reforme}
 \left\{
\begin{aligned}
  &\eta_r+\dfrac{R\theta_+r^{n-1}}{\epsilon \mu v_+^2}\eta = \dfrac{Rr^{n-1}}{\epsilon  \mu v_+}\chi +F(\eta,\chi),
\\[2mm]
&\chi_r=-\frac{\alpha}{\kappa r^{n-1}}+ 
\frac{\epsilon c_P}{\kappa r^{n-1}}\chi + 
\frac{\epsilon^3}{2}\frac{(v_++\eta)^2}{\kappa r^{3(n-1)}}\\
&\hspace{3cm}-\epsilon^2\mu \dfrac{(v_++\eta)\eta_r}{\kappa r^{2(n-1)}}
+\frac{1}{\kappa r^{n-1}}\int^\infty_r\Phi(\eta,\eta_r)(s)\,ds,
\quad r>1,\\[2mm]
&\displaystyle{\lim_{r\to \infty}}\eta(r) =0,\quad \displaystyle{\lim_{r\to \infty}}\chi(r) =0,  \\[2mm]
  &\eta(1)= \eta_-:=v_--v_+,\quad \chi(1)= \chi_-:=\theta_--\theta_+,\quad  (u_->0),\\
&\chi(1)= \chi_-,\quad (u_- < 0). 
 \end{aligned}
 \right.
\end{align}
Here, note that the integral constant $\alpha$ is uniquely determined by
$\eta$ and $\chi$ later, by a formula in (3.1) for the inflow problem, and 
 in (4.3) for the outflow problem, so that it holds $\chi(1)=\chi_-$. 
Once the desired solution $(\eta, \chi)$ of (\ref{reforme}) is obtained, 
the velocity $u$ is immediately obtained by (\ref{utov}) and (\ref{transec}) as
\begin{equation}
u(r)=\frac{u_-(v_++\eta(r))}{v_-r^{n-1}}\quad (u_->0),
\quad
u(r)=\frac{u_-(v_++\eta(r))}{(v_++\eta(1))r^{n-1}}\quad (u_- < 0).
\end{equation} 
The theorem for the reformulated problem (\ref{reforme}) which we need to prove is
\noindent
\begin{thm} 
\label{rmt}
Let $n \ge 3$. Then, for any $v_+>0$ and $\theta_+>0$, there exist positive constants $\epsilon_0$ and $C$ satisfying the followings:

\medskip

\noindent
{\rm (I)} \ Let $u_->0$. If $|u_-|+|\eta_-|+|\chi_-|\le \epsilon_0$, 
there exists a unique smooth solution $(\eta,\chi)$ of the problem {\rm (\ref{reforme})} satisfying 
\begin{align}
\label{decayi}
|\eta(r)|,\ \ |\chi(r)|\le Cr^{-(n-2)}(|u_-|+|\eta_-|+|\chi_-|),\qquad  r\ge 1. 
\end{align}
{\rm (II)} \ Let $u_- < 0$. If $|u_-|+|\chi_-| \le \epsilon_0$, 
there exists a unique smooth solution $(\eta,\chi)$ of the problem {\rm (\ref{reforme})} satisfying
\begin{align}
\label{decayo}
|\eta(r)|,\ \ |\chi(r)|\le Cr^{-(n-2)}(|u_-|+|\chi_-|), \qquad r\ge 1. 
\end{align}
\end{thm}
\noindent
We can easily see that the main Theorem \ref{mt} is a direct  consequence of Theorem \ref{rmt}.

\section{Inflow problem}
In this section, 
we consider the case $u_->0$, that is, the inflow problem, and
show the result (I) in Theorem \ref{rmt}. 
In this case, note $\epsilon = u_-/v_->0$.  
Solving the differential equations (\ref{reforme}) in terms of $\eta$ and $\chi$ 
with the boundary conditions $\eta(1)=\eta_-$ and $\chi(1)=\chi_-$,  
provided the right hand sides of the equations are known,
we have the following differential-integral equations:
\begin{equation}
\label{refint}
\left\{\begin{array}{l}
\eta(r)=G(r,1)\eta_-
+\dfrac{R}{\epsilon \mu v_+}\int_{1}^rG(r,s)\chi(s)s^{n-1}\,ds
+\int_{1}^rG(r,s)F(\eta,\chi)(s)\,ds,
\\[10pt]
\chi(r) = \dfrac{\alpha}{\kappa(n-2)r^{n-2}}+H(\eta,\eta_r,\chi)(r),
\\[10pt]
\alpha = \kappa(n-2)\chi_--\kappa(n-2)
H(\eta,\eta_r,\chi)(1),\\[10pt]
(\eta_r)(r) =-\dfrac{R\theta_+r^{n-1}}{\epsilon\mu v_+^2}\eta(r)+
\dfrac{Rr^{n-1}}{\epsilon \mu v_+}\chi(r) +F(\eta,\chi)(r),\quad r \ge 1,
\end{array}
 \right.\,
\end{equation}
where 
\begin{align}
\label{GKH}
\begin{aligned}
&G(r,s) := e^{-\frac{\omega}{\mu\epsilon}(r^n-s^{n})},\qquad \omega:=\frac{R\theta_+}{v_+^2n},\\
&H(\eta,\eta_r,\chi):= -\int_r^\infty \dfrac{1}{\kappa \tau^{n-1}}
\Big\{{\epsilon c_P}\chi + 
\frac{\epsilon^3}{2}\frac{(v_++\eta)^2}{\tau^{2(n-1)}}\\
&\hspace{5cm} -\epsilon^2\mu \dfrac{(v_++\eta)\eta_r}{\tau^{(n-1)}}
+\int^\infty_\tau\Phi(\eta,\eta_r)(s)\,ds\Big\}\,d\tau.
\end{aligned}
\end{align}
Here, note that we can easily derive the third equality for $\alpha$ in (\ref{refint}) by substituting $r=1$ into the second equation of (\ref{refint}). 
To prove the existence of the solution of (\ref{refint}), we introduce the Banach space $X_l$ with its norm $\|\cdot\|_{X_l}$ as
\begin{align}
\label{bsn}
X_l=\{\eta \in C([1,\infty));\ \sup_{r\ge 1}|r^l\eta(r)|< \infty \},
\quad \|\eta\|_{X_l} = \sup_{r\ge 1}|r^l\eta(r)|.
\end{align} 
To find the solution of (\ref{reforme}) with the decay rate estimate (\ref{decayi}), we look for a solution of (\ref{refint}) in the Banach space $X_{n-2}$, 
with its norm $\|\cdot\|_{X_{n-2}}$.
To do that, we construct the approximate sequence $\{\eta^{(m)}\}_{m\ge 0}$, $\{\chi^{(m)}\}_{m\ge 0}$ and $\{\zeta^{(m)}\}_{m\ge 1}$ by
\begin{equation}
\label{sequence}
\left\{\begin{array}{l}
\eta^{(0)}(r)=G(r,1)\eta_-,\qquad \chi^{(0)}(r)=\dfrac{\chi_-}{r^{n-2}}, \\
\eta^{(m+1)}(r)=G(r,1)\eta_-
+\dfrac{R}{\epsilon  \mu v_+}\int_{1}^rG(r,s)\chi^{(m+1)}(s)s^{n-1}\,ds
\\[10pt]
\hspace{6.5cm}+\int_{1}^rG(r,s)F(\eta^{(m)},\chi^{(m)})(s)\,ds,
\\[10pt]
\chi^{(m+1)}(r) = \dfrac{\chi_--H(\eta^{(m)},\zeta^{(m+1)},\chi^{(m)})(1)}{r^{n-2}}+H(\eta^{(m)}, \zeta^{(m+1)},\chi^{(m)})(r),
\\[10pt]
\zeta^{(m+1)}(r) =-\dfrac{R\theta_+r^{n-1}}{\epsilon \mu v_+^2}\eta^{(m)}+ 
\dfrac{Rr^{n-1}}{\epsilon  \mu v_+}\chi^{(m)} +F(\eta^{(m)},\chi^{(m)})(r),
\end{array}
 \right.\,
\end{equation}
where we introduce a new symbol $\zeta^{(m)}(r):=(\eta_r)^{(m)}(r),\,m\ge 1$. 
Here, we note that all $(\eta^{(m)}, \chi^{(m)}), m\ge 0,$ satisfies the desired boundary and 
far field conditions, provided $(\eta^{(m)}, \chi^{(m)}) \in X_{n-2}, m\ge 0$. 
Then, to show $\{\eta^{(m)}\}_{m\ge 0}$ and $\{\chi^{(m)}\}_{m\ge 0}$  are  Cauchy sequences in $X_{n-2}$ for suitably small $|u_-|+|\eta_-|+|\chi_-|$, we prepare
the following lemma. 
\begin{lem}
\label{hyoka}
{\rm (I)} (refer to \cite{HM1}) If $\epsilon\le \frac{n\omega}{(n-2)\mu}$, then it holds that  
\begin{align}
\label{hyoka1}
r^{n-2}G(r,1)=r^{n-2}e^{-\frac{\omega}{\mu\epsilon}(r^n-1)}\le 1, \qquad r\ge 1.
\end{align}
{\rm (II)}
If  $l\ge -n$, 
then there exists a positive constant $C_0$ which is independent of $\epsilon$ and $\mu$ such that for $f \in X_l$,
\begin{align}
\label{ineqG}
\begin{aligned}
&r^{l+n-1}|\int_1^r G(r,s) f(s) ds|
\le C_0\epsilon \mu\|f\|_{X_l},\quad r\ge 1.
\end{aligned}
\end{align}
\end{lem}

\medskip

\noindent
{\it Proof of {\rm (II)}.} \quad  
By noting that  
\begin{align}
e^{-\frac{\omega}{\epsilon\mu}(r^n-s^n)}=(e^{-\frac{\omega}{\epsilon \mu}(r^n-s^n)})_s\cdot \frac{\epsilon \mu}{n\omega s^{n-1}},
\end{align}
we have
\begin{align}
\label{hyo1}
\begin{aligned}
&|\int_1^r e^{-\frac{\omega}{\epsilon \mu}(r^n-s^n)} f(s) ds|\\
&=|\int_1^r (e^{-\frac{\omega}{\epsilon \mu}(r^n-s^n)})_s\cdot \frac{\epsilon \mu}{n\omega s^{n-1}}\cdot s^{-l}\cdot s^{l}\cdot f(s)ds|\\
&\le \frac{\epsilon \mu}{n\omega}\int_1^r (e^{-\frac{\omega}{\epsilon \mu}(r^n-s^n)})_s\cdot s^{-l-(n-1)} ds \cdot \|f\|_{X_l} \\
&=\frac{\epsilon \mu}{n\omega}\{[e^{-\frac{\omega}{\epsilon \mu}(r^n-s^n)}s^{-l-(n-1)}]_1^r\\
&\qquad -\int_1^r e^{-\frac{\omega}{\epsilon \mu}(r^n-s^n)}(-l-(n-1))\cdot s^{-l-n} ds\}\cdot \|f\|_{X_l}\\
&\le \frac{\epsilon \mu}{n\omega}\{\frac{1}{r^{l+n-1}}+(l+n-1)\int_1^r e^{-\frac{\omega}{\epsilon \mu}(r^n-s^n)}\cdot s^{-l-n} ds\}\cdot \|f\|_{X_l}.
\end{aligned}
\end{align}
Here, we divide the integral domain $[1, r]$ into two parts as 
\begin{align}
\label{I1I2}
\int_1^r e^{-\frac{\omega}{\epsilon \mu}(r^n-s^n)}\cdot s^{-l-n} ds=\int_1^{\frac{1+r}{2}} \ \  ds+\int_{\frac{1+r}{2}}^r \ \  ds:= I_1+I_2.
\end{align}
From the assumption $l+n\ge 0$, we see that 
\begin{align}
\label{hyo2}
\begin{aligned}
I_1&=\int_1^{\frac{1+r}{2}} e^{-\frac{\omega}{\epsilon \mu}(r^n-s^n)}\cdot \frac{1}{s^{l+n}} ds\le \int_1^{\frac{1+r}{2}} e^{-\frac{\omega}{\epsilon \mu}(r^n-s^n)}\cdot 1ds \\
&\le \int_1^{\frac{1+r}{2}} e^{-\frac{n\omega}{2\epsilon \mu}(r-s)}\cdot e^{-\frac{n\omega}{2\epsilon \mu}(r-s)}ds
\le e^{-\frac{n\omega}{4\epsilon \mu}(r-1)}\cdot \int_1^{\frac{1+r}{2}} e^{-\frac{n\omega}{2\epsilon \mu}(r-s)}ds\\
&\le e^{-\frac{n\omega}{4\epsilon \mu}(r-1)}\cdot \frac{2\epsilon \mu}{n\omega}\le \frac{C_0\mu \epsilon}{r^{l+n}},\\
I_2&=\int_{\frac{1+r}{2}}^r e^{-\frac{\omega}{\epsilon \mu}(r^n-s^n)}\cdot \frac{1}{s^{l+n}} ds \le (\frac{2}{1+r})^{l+n}\int_{\frac{1+r}{2}}^r e^{-\frac{\omega}{\epsilon \mu}(r^n-s^n)}ds \\
&\le (\frac{2}{1+r})^{l+n}\frac{\epsilon\mu}{n\omega}\le\frac{ C_0\mu\epsilon}{r^{l+n}}.
\end{aligned}
\end{align}
Substituting (\ref{I1I2}) and (\ref{hyo2}) into (\ref{hyo1}), we obtain the desired estimate
\begin{align*}
|\int_1^r e^{-\frac{\omega}{{\mu}\epsilon}(r^n-s^{n})} f(s) ds|\le C_0\epsilon \mu r^{-l-(n-1)}\cdot \|f\|_{X_l}.
\end{align*}
Hence, the proof of Lemma \ref{hyoka} is completed. \qquad $\Box$ 

\medskip

To estimate $\{\eta^{(m)}\}_{m\ge 0}$ and $\{\chi^{(m)}\}_{m\ge 0}$ smoothly, we also prepare the following proposition. 
\begin{pro}
\label{subhyoka}

Assume that 
\begin{align}
\label{Mbound}
\|\eta^{(m)}\|_{X_{n-2}}\le M,\quad \|\chi^{(m)}\|_{X_{n-2}}\le M,\quad M\le \frac{v_+}{2},\quad 0<u_- \le 1, 
\end{align}
for a positive constant $M\le 1$. Then, there exists a positive constant 
$C_*$ which depends on $v_+$ and $\theta_+$, but not on $M, u_-, \epsilon, \mu$, and  $\kappa$, such that the following estimates hold: 

\medskip

\noindent
{\rm (I)} $|F(\eta^{(m)},\chi^{(m)})(r)| 
\le \dfrac{C_*(M^2+u_-)}{\epsilon \mu \ r^{n-3}}$, \qquad 
{\rm (II)} $|\zeta^{(m+1)}(r)| \le  \dfrac{C_*(M+u_-)}{\epsilon \mu r^{-1}}$,

\medskip

\noindent
{\rm (III)} $|\Phi(\eta^{(m)},\zeta^{(m+1)})(r)|\le \dfrac{C_* u_- (M+(1+\mu) u_-)}{r^{n-1}}$,
 
\medskip
 
\noindent
{\rm (IV)} $|H(\eta^{(m)},\zeta^{(m+1)},\chi^{(m)})(r)|\le \dfrac{C_*(1+\mu+M)u_-}
{\kappa r^{n-2}}$,
\quad $r \ge 1$.
\end{pro}

\noindent
\proof  

\noindent
\underline{Estimate of (I).}\quad
By the definition of $F$ in (\ref{etaF}), we obtain that 
\begin{align}
\label{estF}
\begin{aligned}
|F(\eta^{(m)},\chi^{(m)})(r)|&\le \frac{C_*M^2}{\epsilon {\mu}r^{n-3}}+\dfrac{C_*\epsilon}{\mu r^{n-1}}+\dfrac{MC_*\epsilon}{\mu r^{2n-3}}\le \dfrac{C_*(M^2+u_-^2M+u_-^2)}{\epsilon \mu r^{n-3}}\\
&\le \dfrac{C_*(M^2+u_-)}{\epsilon \mu  r^{n-3}}.
\end{aligned}
\end{align}

\noindent
\underline{Estimate of (II).}\quad
From the definition of $\zeta^{(m+1)}(r)$ in (\ref{sequence}), the assumption (\ref{Mbound}), and the estimate (\ref{estF}), 
we have 
\begin{align}
\label{zetam}
\begin{aligned}
|\zeta^{(m+1)}(r)| \le \dfrac{C_*Mr}{\epsilon \mu}+|F(\eta^{(m)},\chi^{(m)})(r)|&\le \dfrac{C_*(Mr+M^2+u_-)}{\epsilon \mu}\\
&\le \dfrac{C_*(M+u_-)}{\epsilon \mu r^{-1}}.
\end{aligned}
\end{align}

\noindent
\underline{Estimate of (III).}\quad
From the definition of (\ref{Phi}), we see that 
\begin{align}
\begin{aligned}
&\Phi(\eta^{(m)},\zeta^{(m+1)})(r)\\
&=-4\epsilon^2\nu (n-1)\dfrac{(v_++\eta^{(m)})\zeta^{(m+1)}}{r^{n}}+2\epsilon^2\nu n(n-1)\dfrac{(v_++\eta^{(m)})^2}{r^{n+1}}. 
\end{aligned}
\end{align}
Applying the estimate of $\zeta^{(m+1)}(r)$ in (\ref{zetam}) and noting (\ref{mutonu}), we obtain
\begin{align}
\label{Phim}
\begin{aligned}
&|\Phi(\eta^{(m)},\zeta^{(m+1)})(r)| 
\le \dfrac{C_*\epsilon\nu(M+u_-)}{r^{n-1} \mu}
+\dfrac{C_*\epsilon^2\nu}{r^{n+1}}
\le \frac{C_*u_-(M+(1+\mu)u_-)}{r^{n-1}}.
\end{aligned}
\end{align}

\bigskip

\noindent
\underline{Estimate of (IV).}\quad By the definition of $H(\eta,\eta_r,\chi)$ in (\ref{GKH}), we see that 
\begin{align}
\label{h-m}
\begin{aligned}
&H(\eta^{(m)},\zeta^{(m+1)},\chi^{(m)})(r)\\
&= -\int_r^\infty \dfrac{1}{\kappa \tau^{n-1}}
\Big\{{\epsilon c_P}\chi^{(m)} + 
\frac{\epsilon^3}{2}\frac{(v_++\eta^{(m)})^2}{\tau^{2(n-1)}}
-\epsilon^2\mu \dfrac{(v_++\eta^{(m)})\zeta^{(m+1)}}{\tau^{(n-1)}}\\
&\qquad +\int^\infty_\tau\Phi(\eta^{(m)},\zeta^{(m+1)})(s)\,ds\Big\}\,d\tau. 
\end{aligned}
\end{align}
Applying the estimates (\ref{zetam}) and (\ref{Phim}), we have
\begin{align}
\label{H1}
\begin{aligned}
&|H(\eta^{(m)},\zeta^{(m+1)},\chi^{(m)})(r)|\\
&\le \int_r^\infty \dfrac{1}{\kappa \tau^{n-1}}
\Big\{{\epsilon c_P}|\chi^{(m)}| + 
\frac{\epsilon^3}{2}\frac{(v_++|\eta^{(m)}|)^2}{\tau^{2(n-1)}}\\
&\qquad +\epsilon^2\mu \dfrac{(v_++|\eta^{(m)}|)|\zeta^{(m+1)}|}{\tau^{(n-1)}}+|\int^\infty_\tau\Phi(\eta^{(m)},\zeta^{(m+1)})(s)\,ds|\Big\}\,d\tau\\
&\le \dfrac{C_*}{\kappa}\int_r^\infty  
\Big\{\frac{|u_-|}{\tau^{n-1}}+\frac{|u_-|^3}{\tau^{3(n-1)}}+\frac{\mu|u_-|^2}{\tau^{2(n-1)}}\cdot\dfrac{1}{\epsilon \mu \tau^{-1}}
+\frac{(M+(1+\mu)u_-) |u_-|}{\tau^{2n-3}} \Big\}\,d\tau \\
&\le \frac{C_*|u_-|(1+\mu+M)}{\kappa r^{n-2}}.
\end{aligned}
\end{align}
From the above estimate (\ref{H1}), we specifically have 
\begin{align}
\label{H2}
\begin{aligned}
|H(\eta^{(k)},\zeta^{(k+1)},\chi^{(k)})(1)|\le C_*\kappa^{-1}(1+\mu+M)|u_-|.
\end{aligned}
\end{align}
Thus, the proof of Proposition \ref{subhyoka} is completed.\hspace{2cm}$\Box$

\medskip

Now, by using Lemma \ref{hyoka}, we show the uniform boundedness of 
$\eta^{(m)}$\ and  \ $\chi^{(m)}\ \ (m\ge 0)$ in $X_{n-2}$ for suitably small 
$|u_-|+|\eta_-|+|\chi_-|$. 
More precisely, we show that for any fixed $v_+$ and $\theta_+$, 
there exist positive constants $\epsilon_0$ and $C$ which are
independent of $u_-$, $\eta_-$ and $\chi_-$ such that if 
$|u_-|+|\eta_-|+|\chi_-|\le \epsilon_0$, then there exists a positive
constant $M$ satisfying
\begin{equation}
\label{ubM}
\|\eta^{(m)}\|_{X_{n-2}}, \quad \|\chi^{(m)}\|_{X_{n-2}} \le M \le C(|u_-|+|\eta_-|+|\chi_-|),\quad m\ge 0.
\end{equation}
In what follows, the letters $C$ and $\epsilon_0$ are used to represent generic
positive constants which are independent of $u_-$, $\eta_-$ and $\chi_-$,
but may depend on $v_+$, $\theta_+$ and other fixed constants like 
$\mu, \nu, \kappa, n,\dots,etc$. In particular, $C_*$ denote constants $C$ 
which are independent of $\mu, \nu$,  and $\kappa$.
For the proof, we first assume
\begin{equation}
\label{sa1}
|\eta_-|=|v_--v_+|\le \frac{v_+}{2}, \quad u_-\le \frac{n\omega v_+}{2(n-2)\mu}, 
\quad u_- \le 1,
\end{equation}
which assures 
\begin{equation}
\label{asve}
\frac{v_+}{2} \le v_- \le \frac{3 v_+}{2},\quad \epsilon=\frac{u_-}{v_-}\le \frac{n\omega}{(n-2)\mu}.
\end{equation}
Let us show (\ref{ubM}) by mathematical induction:

\smallskip

\noindent
\underline{Case $m=0$}.\quad
Because of (\ref{asve}), we have from Lemma \ref{hyoka}-(I) that 
\begin{align}
\label{m0}
\begin{aligned}
&|r^{n-2}\eta^{(0)}(r)|=|G(r,1)\eta_-|\le |\eta_-|,\\
& |r^{n-2}\chi^{(0)}(r)|=|\chi_-|,\quad r\ge 1.
\end{aligned}
\end{align}
Hence, we ask the constant 
$M$ to satisfy 
\begin{align}
\label{sa2}
|\eta_-|,\quad |\chi_-| \le M,
\end{align}
so that it holds $\|\eta^{(0)}\|_{X_{n-2}},\ \|\chi^{(0)}\|_{X_{n-2}} \le M$.

\smallskip

\noindent
\underline{Case $m=k+1$}\ $(k\ge 0)$.\quad
Suppose $\|\eta^{(k)}\|_{X_{n-2}},\ \|\chi^{(k)}\|_{X_{n-2}} \le M$. 
Here, we ask the constant $M$ to satisfy another assumption
\begin{align}
\label{sa3}
M\le \frac{v_+}{2},\quad M\le 1,
\end{align}
which, in particular, imply
\begin{align}
\label{bound}
|\eta^{(k)}(r)|\le \frac{v_+}{2},\quad r \ge 1.
\end{align}
Then, by using Proposition \ref{subhyoka}, we estimate $\chi^{(k+1)}$, $\eta^{(k+1)}$, and $\zeta^{(k+1)}$ defined by (\ref{sequence}). For the estimate of $\chi^{(k+1)}$, by using Proposition \ref{subhyoka}-(IV), we obtain
\begin{align}
\label{chiestimate}
\begin{aligned}
&|r^{n-2}\chi^{(k+1)}| \\
&\le |\chi_--H(\eta^{(k-1)},\zeta^{(k)},\chi^{(k-1)})(1)|+|r^{n-2}H(\eta^{(k-1)},\zeta^{(k)},\chi^{(k-1)})(r)|\\
&\le |\chi_-|+C_*(1+\mu)\kappa^{-1}|u_-|.
\end{aligned}
\end{align}
By using the Lemma \ref{hyoka}-(I), (\ref{chiestimate}), and applying   Proposition \ref{subhyoka}-(I), we can estimate $\eta^{(k+1)}$ as 
\begin{align}
\label{etaestimate}
\begin{aligned}
&|r^{n-2}\eta^{(k+1)}|\le |r^{n-2}G(r,1)\eta_-|
+|\dfrac{Rr^{n-2}}{\epsilon  \mu v_+}\int_{1}^rG(r,s)\chi^{(k+1)}(s)s^{n-1}\,ds|\\
&\hspace{6cm}+|r^{n-2}\int_{1}^rG(r,s)F(\eta^{(k)},\chi^{(k)})(s)|\\
&\le |\eta_-|+\dfrac{C_*r^{n-2}(|\chi_-|+(1+\mu)\kappa^{-1})|u_-|}
{\epsilon  \mu}\int_{1}^rG(r,s)\cdot s\,ds\\
&\hspace{4cm} +\dfrac{C_*(M^2+|u_-|)}{\epsilon \mu}|r^{n-2}\int_{1}^rG(r,s)\frac{1}{s^{n-3}}ds|\\
&\le |\eta_-|+\dfrac{C_*r^{n-2}(|\chi_-|+(1+\mu)\kappa^{-1})|u_-|}
{\epsilon  \mu}\cdot C_0\epsilon \mu r^{-(n-2)}\\
&\hspace{6cm}+\dfrac{C_*r^{n-2}(M^2+|u_-|)}{\epsilon \mu}
\cdot C_0\frac{\epsilon \mu}{r^{2n-4}}\\
&\le |\eta_-|+C_*(|\chi_-|+(1+\mu)\kappa^{-1})|u_-|+C_*(M^2+|u_-|)\\
&\le |\eta_-|+C_*(|\chi_-|+(1+(1+\mu)\kappa^{-1})|u_-|+M^2).
\end{aligned}
\end{align}
Combining the above inequalities (\ref{chiestimate}) and (\ref{etaestimate}), we obtain the following inequality
\begin{align}
\label{hyouka1}
|r^{n-2}\eta^{(k+1)}|,\quad |r^{n-2}\chi^{(k+1)}|\le 
|\eta_-|+C_*(|\chi_-|+(1+(1+\mu)\kappa^{-1})|u_-|+M^2).
\end{align}
\noindent
Therefore, we further assume 
\begin{align}
\label{sa4}
|\eta_-|+C_*(|\chi_-|+(1+(1+\mu)\kappa^{-1})|u_-|) \le \frac{M}{2}, 
\quad  C_*M\le \frac{1}{2},
\end{align}
so that  (\ref{hyouka1}) gives the desired estimate $\|\chi^{(k+1)}\|_{X_{n-2}},\ \|\eta^{(k+1)}\|_{X_{n-2}} \le M$. We can easily see  that there exists a
positive constant $\epsilon_0$ such that if 
$|u_-|+|\eta_-|+|\chi_-| \le \epsilon_0$, all the assumptions 
(\ref{sa1}), (\ref{sa2}), (\ref{sa3}) and (\ref{sa4}) hold. In particular,
$M$ can be chosen by
\begin{align}
\label{M1}
M =2 |\eta_-|+2C_*(|\chi_-|+(1+(1+\mu)\kappa^{-1})|u_-|)
\le C(|\eta_-|+|\chi_-|+|u_-|),
\end{align}
which proves the uniform boundedness of 
$\chi^{(m)}$ and $\eta^{(m)}\ (m\ge 0)$ in $X_{n-2}$ with the estimate (\ref{ubM}). 
Here, we note that if we assume $\mu \in (0,1]$, then 
we can easily see $\epsilon_0$ and $C$ can be chosen independent of $\mu$
by (\ref{sa1}) and (\ref{M1}). 
Finally, the proof to show $\{(\chi^{(m)}, \eta^{(m)})\}_{m\ge 0}$
is a Cauchy sequence in $X_{n-2}\times X_{n-2}$ is very standard. In fact, by using the definition of $\chi^{(m)}$ in (\ref{sequence}), we obtain the following estimate 
\begin{align}
\label{ccc}
\begin{aligned}
&\|\chi^{(m+1)}-\chi^{(m)}\|_{X_{n-2}} \\
&\le\sup_{r\ge 1}|r^{n-2}(H(\eta^{(m)},\zeta^{(m+1)},\chi^{(m)})-H(\eta^{(m-1)},\zeta^{(m)},\chi^{(m-1)}))(r)|\\[2mm]
&\le C|u_-|\left(\|\eta^{(m)}-\eta^{(m-1)}\|_{X_{n-2}}+\|\chi^{(m)}-\chi^{(m-1)}\|_{X_{n-2}}\right), \quad m \ge 1,
\end{aligned}
\end{align}
where we use the following estimate
\begin{align}
\label{zatacauchy1}
\|\zeta^{(m+1)}-\zeta^{(m)}\|_{X_{-1}}\le \frac{C}{\epsilon \mu}(\|\eta^{(m)}-\eta^{(m-1)}\|_{X_{n-2}}+\|\chi^{(m)}-\chi^{(m-1)}\|_{X_{n-2}}).
\end{align}
By using the above estimate (\ref{ccc}) and Lemma \ref{hyoka}, we have the following estimate for $m \ge 1$: 
\begin{align}
\label{cce}
\begin{aligned}
&\|\eta^{(m+1)}-\eta^{(m)}\|_{X_{n-2}} \\
&=\sup_{r\ge 1}|r^{n-2}(\frac{R}{\epsilon\tilde\mu v_+}\int_{1}^rG(r,s)(\chi^{(m+1)}-\chi^{(m)})s^{n-1} \,ds)|\\[2mm]
&\hspace{2cm}+\sup_{r\ge 1}|r^{n-2}(\int_{1}^rG(r,s)(F(\eta^{(m)},\chi^{(m)})-F(\eta^{(m-1)},\chi^{(m-1)})) \,ds)|\\[2mm]
&\le C \|\chi^{(m+1)}-\chi^{(m)}\|_{X_{n-2}} \\
&\quad+C(|\eta_-|+ |u_-|+|\chi_-|)\left(\|\eta^{(m)}-\eta^{(m-1)}\|_{X_{n-2}}+\|\chi^{(m)}-\chi^{(m-1)}\|_{X_{n-2}}\right)\\
&\le C(|\eta_-|+ |u_-|+|\chi_-|)\left(\|\eta^{(m)}-\eta^{(m-1)}\|_{X_{n-2}}+\|\chi^{(m)}-\chi^{(m-1)}\|_{X_{n-2}}\right).
\end{aligned}
\end{align}
Combining the above two estimates (\ref{ccc}) and (\ref{cce}), and taking $\epsilon_0$ suitably small again if needed, 
we can show that 
\begin{align*}
\begin{aligned}
&\left(\|\eta^{(m+1)}-\eta^{(m)}\|_{X_{n-2}}+\|\chi^{(m+1)}-\chi^{(m)}\|_{X_{n-2}}\right)\\
&\hspace{4cm}\le \frac{1}{2}\left(\|\eta^{(m)}-\eta^{(m-1)}\|_{X_{n-2}}+\|\chi^{(m)}-\chi^{(m-1)}\|_{X_{n-2}}\right), \quad m \ge 1,
\end{aligned}
\end{align*}
which proves that $\{(\chi^{(m)}, \eta^{(m)})\}_{m\ge 0}$
is a Cauchy sequence in ${X_{n-2}}\times {X_{n-2}}$, that is, there exist functions $(\eta, \chi)\in {X_{n-2}}\times {X_{n-2}}$ such that $\eta^{(m)}\to \eta$, and $\chi^{(m)}\to \chi$ in $X_{n-2}$. In addition, we can see from (\ref{zatacauchy1}) that $\{\zeta^{(m)}\}_{m\ge 1}$ is also a convergent sequence in $X_{-1}$, which implies that there exists a function $\zeta \in {X_{-1}}$ such that $\zeta^{(m)}\to \zeta$ in $X_{-1}$. 
Therefore, we can see from (\ref{sequence}) that the limit functions $(\eta, \chi, \zeta)$ satisfy the following equations
\begin{equation}
\label{sequencelimit}
\left\{\begin{array}{l}
\eta(r)=G(r,1)\eta_-
+\dfrac{R}{\epsilon  \mu v_+}\int_{1}^rG(r,s)\chi(s)s^{n-1}\,ds
\\[10pt]
\hspace{6.5cm}+\int_{1}^rG(r,s)F(\eta,\chi)(s)\,ds,
\\[10pt]
\chi(r) = \dfrac{\chi_--H(\eta,\zeta,\chi)(1)}{r^{n-2}}+H(\eta, \zeta,\chi)(r),
\\[10pt]
\zeta(r) =-\dfrac{R\theta_+r^{n-1}}{\epsilon \mu v_+^2}\eta+ 
\dfrac{Rr^{n-1}}{\epsilon  \mu v_+}\chi +F(\eta,\chi)(r),
\end{array}
 \right.\,
\end{equation}
where we can easily see that $\zeta(r)=(\eta_r)(r)$. Thus, the solution $(\chi, \eta)$ of (\ref{reforme}) with the 
desired estimate (\ref{decayi}) is obtained. 
The arguments on the regularity and uniqueness of the solution 
are standard, so we omit them. $\Box$

\section{Outflow problem}
This subsection is devoted to considering the case $u_-<0$, that is, the outflow problem, and we show the result (II) in Theorem \ref{mt}. 
In this case, recalling 
\begin{align}
\label{epout}
\epsilon := u_-/(v_++\eta(1))<0, 
\end{align}
we rewrite the first equation of $\eta$ in (\ref{reforme}) in the form:
\begin{equation}
\label{etar}
\begin{array}{l}
\eta_r+\dfrac{R\theta_+r^{n-1}}{u_- \mu v_+}\eta =-\dfrac{R\theta_+\eta(1)}{u_-{\mu}v_+^2}r^{n-1}\eta+ \dfrac{Rr^{n-1}}{\epsilon  \mu v_+}\chi +F(\eta,\chi),\\
\end{array}
\end{equation}
where $F(\eta,\chi)$ is defined by the same form as (\ref{etaF}) with the definition (\ref{epout}). 
Noting $u_-<0$ and the far field condition $\eta(\infty)=0$, we can solve equation (\ref{etar}) in terms of $\eta$ by applying Duhamel's principle, provided the right hand side is known. 
On the other hand, the integrated equation of $\chi$ is the same as the second equation of (\ref{refint}) in Section 3, so our purpose in this section is to obtain the solution of the following differential-integral equations (\ref{outequ}) in terms of $\eta$ and $\chi$ with the boundary condition $\chi(1)=\chi_-$:
\begin{equation}
\label{outequ}
\left\{\begin{array}{l}
\eta(r)=-\int_{r}^{\infty} \ti G(r,s)\left(-\dfrac{R\theta_+\eta(1)}{u_-{\mu}v_+^2}s^{n-1}\eta+ \dfrac{Rs^{n-1}}{\epsilon \mu v_+}\chi +F(\eta,\chi)\right)(s)\,ds,
\\[10pt]
\chi(r) = \dfrac{\alpha}{\kappa(n-2)r^{n-2}}+H(\eta,\eta_r,\chi)(r),
\\[10pt]
\alpha = \kappa(n-2)(\theta_--\theta_+)-\kappa(n-2)
H(\eta,\eta_r,\chi)(1),\\[10pt]
(\eta_r)(r) =-\dfrac{R\theta_+r^{n-1}}{\epsilon \mu v_+^2}\eta(r)+
\dfrac{Rr^{n-1}}{\epsilon  \mu v_+}\chi(r) +F(\eta,\chi)(r),\quad r \ge 1,
\end{array}
 \right.\,
\end{equation}
where $H$ is defined in the same form as in (\ref{GKH}) and
\begin{align}
\begin{aligned}
\ti G(r,s) := e^{-\frac{\bar\omega}{|u_-|{\mu}}(s^n-r^{n})},\qquad  \bar\omega:=\frac{R\theta_+}{v_+ n}(>0).
\end{aligned}
\end{align}
To prove the existence of the solution of (\ref{outequ}) 
with the decay rate estimate (\ref{decayo}),
we look for a solution of (\ref{outequ}) in the Banach space $X_{n-2}$ defined in (\ref{bsn}) with $l=n-2$. In the same way as in the last section, 
we construct the approximate sequence $\{\eta^{(m)}\}_{m\ge 0}$, $\{\chi^{(m)}\}_{m\ge 0}$ and $\{\zeta^{(m)}\}_{m\ge 1}$ corresponding to (\ref{outequ}) by
\begin{equation}
\label{sequenceout}
\left\{\begin{array}{l}
\eta^{(0)}:=0, \qquad \chi^{(0)}:=0,\\[10pt]
\eta^{(m+1)}(r)=-\int_{r}^{\infty} \ti G(r,s)\Big(-\dfrac{R\theta_+\eta^{(m)}(1)}{u_-{\mu}v_+^2}s^{n-1}\eta^{(m)}+ \dfrac{Rs^{n-1}}{\epsilon \mu v_+}\chi^{(m+1)} \\[10pt]
\hspace{8.5cm}+ F(\eta^{(m)},\chi^{(m)})\Big)(s)\,ds,\\[10pt]
\chi^{(m+1)}(r) = \dfrac{\chi_--H(\eta^{(m)},\zeta^{(m+1)},\chi^{(m)})(1)}{r^{n-2}}+H(\eta^{(m)}, \zeta^{(m+1)},\chi^{(m)})(r),
\\[10pt]
\zeta^{(m+1)}(r) =-\dfrac{R\theta_+r^{n-1}}{\epsilon \mu v_+^2}\eta^{(m)}+
\dfrac{Rr^{n-1}}{\epsilon  \mu v_+}\chi^{(m)} +F(\eta^{(m)},\chi^{(m)})(r),
\end{array}
 \right.\,
\end{equation}
where $\zeta^{(m)}(r):=(\eta_r)^{(m)}(r),\,m\ge 1$. 
Here, we note that although $\chi^{(0)}$ doesn't satisfy the boundary condition 
$\chi^{(0)}(1)= \chi_-$, 
all $(\eta^{(m)}, \chi^{(m)}), m\ge 1,$ satisfy the desired boundary and 
far field conditions, provided $(\eta^{(m)}, \chi^{(m)}) \in X_{n-2}, m\ge 1.$ 
Then, to show $\{\eta^{(m)}\}_{m\ge 0}$ and $\{\chi^{(m)}\}_{m\ge 0}$ form Cauchy sequences in $X_{n-2}$ for 
suitably small $|u_-|$ and $|\chi_-|$, we prepare the following lemma.
\begin{lem}
\label{sekibunhyokaout}
If  $l\ge -n$, then there exists a positive constant $C_0$ which is independent of $\epsilon$ and $\mu$ such that  for $g \in X_{l}$,
\begin{align*}
\label{hyoka2''}
r^{l+n-1}|\int_{r}^{\infty}\ti G(r,s)g(s)\,ds|
\le C_0\mu|u_-|\cdot \sup_{s \ge r}|s^{l}g(s)|,\quad r\ge 1.
\end{align*}
\end{lem}
\proof \ \ Because $g \in X_l$, it holds that
\begin{align*}
\begin{aligned}
&|\int_{r}^{\infty}\ti G(r,s)g(s)\,ds|\,  \le\int_r^{\infty} |e^{-\frac{\bar\omega}{|u_-| \mu}(s^n-r^n)} g(s)| ds \\
&\le \int_r^{\infty} e^{-\frac{\bar\omega}{|u_-| \mu}(s^n-r^n)} s^{-l} ds \cdot \sup_{s \ge r}|s^{l}\cdot g(s)|\\
=-&\frac{|u_-| \mu}{n\bar\omega}\int_r^{\infty} (e^{-\frac{\bar\omega}{|u_-| \mu}(s^n-r^n)})_s\cdot s^{-l-(n-1)} ds\cdot \sup_{s\ge r}|s^l\cdot g(s)|\\
=-&\frac{|u_-| \mu}{n\bar\omega}\Big\{[e^{-\frac{\bar\omega}{|u_-| \mu}(s^n-r^n)}\cdot s^{-l-(n-1)}]_r^{\infty}\\
&\quad-\int_r^{\infty} e^{-\frac{\bar\omega}{|u_-| \mu}(s^n-r^n)}(-l-(n-1))\cdot s^{-l-n} ds\Big\}\cdot \sup_{s\ge r}|s^l\cdot g(s)| \\
=& -\frac{|u_-| \mu}{n\bar\omega}\{\frac{-1}{r^{l+n-1}}+(l+n-1)\int_r^{\infty} e^{-\frac{\bar\omega}{|u_-| \mu}(s^n-r^n)}\cdot s^{-l-n} ds\}\cdot \sup_{s\ge r}|s^l\cdot g(s)|\\
\le& \frac{|u_-| \mu}{n\bar\omega}\{\frac{1}{r^{l+n-1}}+\frac{|l+n-1|}{r^{l+n}} \int_r^{\infty} e^{-\frac{\bar\omega}{|u_-| \mu}(s^n-r^n)}\ ds\}\cdot \sup_{s\ge r}|s^l\cdot g(s)|\\
\le& C_0|u_-|\mu \ \frac{1}{r^{l+n-1}} \cdot \sup_{s\ge r}|s^l\cdot g(s)|. 
\end{aligned}
\end{align*}
Thus, the proof of Lemma \ref{sekibunhyokaout} is completed.$\qquad \Box
$


\medskip

By using Lemma \ref{sekibunhyokaout}, we show that for any fixed $v_+$ and $\theta_+$, 
there exist positive constants $\epsilon_0$ and $C$ such that if 
$|u_-|,\ |\chi_-|\le \epsilon_0$, then there exists a positive
constant $M$ satisfying
\begin{equation}
\label{boundout}
\|\eta^{(m)}\|_{X_{n-2}},\quad \|\chi^{(m)}\|_{X_{n-2}} \le M \le C(|u_-|+|\chi_-|),\quad m\ge 0.
\end{equation}
Here and hereafter, the letters $C$ and $\epsilon_0$ are used to represent generic
positive constants which are independent of $u_-$ and $\chi_-$,
but may depend on $v_+$, $\theta_+$ and other fixed constants like 
$\mu, \nu, \kappa, n,\dots,etc$. In particular, $C_*$ denote constants $C$ 
which are independent of $\mu, \nu, \kappa$ as in the previous subsection. 
To prove (\ref{boundout}), we prepare the following proposition.
\begin{pro}
\label{subhyokaout}
Assume that 
\begin{align}
\|\eta^{(m)}\|_{X_{n-2}}\le M,\quad \|\chi^{(m)}\|_{X_{n-2}}\le M,\quad M\le \frac{v_+}{2},\quad |u_-| \le 1,
\end{align}
for a positive constant $M\le 1$. Then, there exists a positive constant
$C_*$ such that the following estimates hold : 

\medskip

\noindent
{\rm (I)}\ $|F(\eta^{(m)},\chi^{(m)})(r)| \le \dfrac{C_*(M^2+|u_-|)}{|u_-| \mu \ r^{n-3}},$ \qquad
{\rm (II)}\ $|\zeta^{(m+1)}(r)|\le \dfrac{C_*(M+|u_-|)}{|u_-|\mu r^{-1}}$,

\medskip

\noindent
{\rm (III)}\  $|\Phi(\eta^{(m)},\zeta^{(m+1)})(r)|\le \dfrac{C_*|u_-|(M+(1+\mu) |u_-|)}{r^{n-1}}$, 

\medskip

\noindent
{\rm (IV)}\  $|H(\eta^{(m)},\zeta^{(m+1)},\chi^{(m)})(r)|\le 
\dfrac{C_*|u_-|(1+\mu + M)}{\kappa r^{n-2}},\quad r \ge 1.$
\end{pro}
Once we note the definition of $\epsilon$ in (\ref{epout}), the proof of Proposition \ref{subhyokaout} is almost the same as that of Proposition \ref{subhyoka}, so we omit the proof. Let us show (\ref{boundout}) by mathematical induction:

\smallskip

\noindent
\underline{Case $m=0$}.\quad
Since $\eta^{(0)}=\chi^{(0)}=0$, then (\ref{boundout}) trivially holds.

\medskip

\noindent
\underline{Case $m=k+1$}\ ($k\ge 0)$.\quad 
Suppose $\|\eta^{(k)}\|_{X_{n-2}}, \quad \|\chi^{(k)}\|_{X_{n-2}}\le M$. 
Here we ask that the constant $M$ satisfies another assumption:
\begin{align}
\label{outuniformk}
M\le \frac{v_+}{2}, \quad M \le 1
\end{align}
which, in particular, implies
\begin{align}
|\eta^{(k)}(1)| \le \|\eta^{(k)}\|_{X_{n-2}}
\le M \le \frac{v_+}{2}.
\end{align}
By applying (IV) in Proposition \ref{subhyokaout}, we estimate $\chi^{(k+1)}$ defined by (\ref{sequenceout})  as 
\begin{align}
\label{uniformkchi}
\begin{aligned}
|r^{n-2}\chi^{(k+1)}| 
&\le |\chi_-|+\sup_{r\ge 1}|r^{(n-2)}H(\eta^{(k-1)},(\eta_r)^{(k)},\chi^{(k-1)})(r)|\\
&\le C_*(|\chi_-|+(1+\mu)\kappa^{-1})|u_-|).
\end{aligned}
\end{align}
Next, we estimate $\eta^{(k+1)}$ defined in (\ref{sequenceout}) as
\begin{align}
\label{outI123}
\begin{aligned}
&|r^{n-2}\eta^{(k+1)}|\\
&\le|r^{n-2}\int_{r}^{\infty} \ti G(r,s)\left(-\dfrac{R\theta_+\eta^{(m)}(1)}{u_-{\mu}v_+^2}s^{n-1}\eta^{(m)}+ \dfrac{Rs^{n-1}}{\epsilon \mu v_+}\chi^{(m+1)} \right..\\
&\hspace{9cm}+F(\eta^{(m)},\chi^{(m)})\Big)(s)\,ds|\\
&=: I_1+I_2+I_3,
\end{aligned}
\end{align}
\begin{align}
\label{outI1}
\begin{aligned}
&I_1:=r^{n-2}\int_{r}^{\infty} \ti G(r,s)|\dfrac{R\theta_+\eta^{(m)}(1)}{u_-{\mu}v_+^2}s^{n-1}\eta^{(m)}|ds \\
&\le C_*r^{n-2}\int_{r}^{\infty} \ti G(r,s)\dfrac{M^2}{|u_-|{\mu}}s\ ds 
\le r^{n-2}\cdot \frac{C_*|u_-| \mu}{ r^{n-2}}\cdot \dfrac{C_0M^2}{|u_-|{\mu}}\le C_*M^2,
\end{aligned}
\end{align}
\begin{align}
\label{outI2}
\begin{aligned}
&I_2:= r^{n-2}\int_{r}^{\infty} \ti G(r,s) |\dfrac{Rs^{n-1}}{\epsilon  \mu v_+}\chi^{(m+1)}| ds\\
&\le r^{n-2}\int_{r}^{\infty} \ti G(r,s) |\dfrac{Rs^{n-1}(v_++\eta(1))}{u_- \mu v_+}|
\cdot\frac{C_*(|\chi_-|+(1+\mu)\kappa^{-1})|u_-|}{s^{n-2}} ds \\
&\le C_*(|\chi_-|+(1+\mu)\kappa^{-1}|u_-|), 
\end{aligned}
\end{align}
\begin{align}
\label{outI3}
\begin{aligned}
&I_3:=  r^{n-2}\int_{r}^{\infty} \ti G(r,s) |F(\eta^{(m)},\chi^{(m)})(s)| ds\\
&\le r^{n-2}\int_{r}^{\infty} \ti G(r,s)\dfrac{C_*(M^2+|u_-|)}{|u_-| \mu \ s^{n-3}} \ ds \\
&\le C_*\frac{(M^2+|u_-|)}{|u_-| \mu} \cdot |u_-|\mu \le C_*(M^2+|u_-|),
\end{aligned} 
\end{align}
where we used Lemma \ref{sekibunhyokaout} and Proposition \ref{subhyokaout}.
Substituting (\ref{outI1})-(\ref{outI3}) into (\ref{outI123}), we obtain
%
%
the estimate of $\chi^{(k+1)}$ and $\eta^{(k+1)}$ as 
\begin{align}
\label{estetaout}
\begin{aligned}
|r^{n-2}\chi^{(k+1)}|, \quad |r^{n-2}\eta^{(k+1)}|\le C_*(|\chi_-|+(1+(1+\mu)\kappa^{-1}) |u_-|+M^2).
\end{aligned}
\end{align}
Therefore, we further assume 
\begin{align}
\label{outuniformall}
C_*(|\chi_-|+(1+(1+\mu)\kappa^{-1}) |u_-|)\le \frac{M}{2},\quad C_*M\le \frac{1}{2}.
\end{align}
so that (\ref{uniformkchi}) and (\ref{estetaout}) give the desired estimate $\|\eta^{(k+1)}\|_{X_{n-2}},\quad  \|\chi^{(k+1)}\|_{X_{n-2}} \le M$.
It is easy to see there exists a positive constant $\epsilon_0$ such that if 
$|u_-|+|\chi_-| \le \epsilon_0$, all the assumptions 
 (\ref{outuniformk})  and (\ref{outuniformall}) hold, and in particular,
$M$ can be chosen by
\begin{align*}
\label{M}
M= 2C_*(|\chi_-|+(1+(1+\mu)\kappa^{-1}) |u_-|) \le C(|\chi_-|+|u_-|),
\end{align*}
which proves the uniform boundedness of 
$\eta^{(m)}$ and $\chi^{(m)}\ (m\ge 0)$ in $X_{n-2}$ with the estimate (\ref{boundout}). 
Once the uniform estimate (\ref{boundout}) is obtained, the remaining arguments on the convergence to the limits $\eta$, $\chi$, and $\eta_r$, and the regularity and uniqueness of
the limit are very standard, as in the last section, so we omit the details. 
Thus, the proof for the result (II) in Theorem 2.1 is completed.


\bigskip

\end{document}